\input amstex 
\documentstyle{amsppt}
\magnification=\magstep1 
\hsize=12.6 cm
\vsize=17.5 cm
\hoffset=2.8pc
\voffset=1.8pc
\redefine \Re{\operatorname{Re\,}}
\redefine \Im{\operatorname{Im\,}}
\define \Res{\operatorname{Res\,}}
\topmatter
\head
{\bf ON INCOMMENSURABILITY OF TWO-DIMENSIONAL STANDARD INTEGER LATTICE 
AND A LATTICE GENERATED BY VECTORS, WHICH HAVE ZETA VALUES AS COORDINATES } 
\endhead 
\head 
{The detailed version, \, \bf Part 1 (INTRODUCTION)} 
\endhead 
\author
Leonid A. Gutnik
\endauthor
\vskip-16pt
\dedicatory
\vskip6pt
{\hfill \sl In memory of Professor Alexander Buchstab} 
\enddedicatory 
\endtopmatter
\eightpoint Abstract. 

This is a revised variant of my preprint:

Max-Plank Institut fuer Mathematik, Bonn, 2001, No 16.
\tenpoint

\rightheadtext{On incommensurability of lattices}

\document 
\baselineskip12.5pt
\TagsOnLeft
\vskip.01pt

\heading{\sl Foreword}\endheading

This is the first part of the detailed version of my work ``On linear 
forms with coefficients in $\Bbb N \zeta (1 + \Bbb N),$'' [11]. The goal of 
this work is to give a complete proof of the following theorem. 
\proclaim{Theorem} Let  
$$
\phi \sb i (x\sb 1, x\sb2) = \sum \limits \sb {k=1} \sp2 (2ik - i - k + 2)
\,\zeta (i + k + 1) \,x\sb k, \qquad i = 1, 2.
$$ 
For any $d \in \Bbb R,$ let $\Vert d \Vert$ be the 
distance between $d$ and $\Bbb Z;$ let $\gamma = 43,464412.$ 
There is a positive constant $c$ such that
$$
\Vert \phi \sb1 (x\sb 1,x\sb2)\Vert + \Vert\phi \sb 2 (x\sb 1,x\sb2)\Vert 
\ge c(\vert x\sb 1\vert +\vert x\sb2\vert )\sp {-\gamma},$$ 
where $x\sb 1 \in \Bbb Z, \; x\sb2\in \Bbb Z, \; 
\vert x\sb 1\vert + \vert x\sb2\vert > 0.$
\endproclaim 

\noindent
{\tt Remark.} The lattices $L_1, \, L_2$ in ${\Bbb R}^n$ are 
said to be incommensurable, if $L_1 \cap L_2 = \{0\}.$ 
The qualitative part of the theorem asserts that lattice generated by the vectors $f_1 = \big(2\,\zeta(3), \, 3\,\zeta(4)\big)$ and $f_2 = 
\big(3\,\zeta(4), \, 6\,\zeta(5)\big)$ and the lattice ${\Bbb Z}^2$ 
(generated by the vectors $(1,0)$ and $(0,1))$ are incommensurable. 
\proclaim{Corollary} \it If $p \in  \Bbb Q, \; q \in  \Bbb Q,$ and 
$p\sp 2 + q\sp 2 > 0,$ then  
$$\{2p\,\zeta (3) + 3q\,\zeta (4), \; 3p\,\zeta (4) + 6q\,\zeta (5)\} 
\nsubseteq \Bbb Q,$$
and therefore  
$$\bigg\{\frac{\zeta(3+2k)}{\zeta(4)}, \; \frac{12\,\zeta(3)\zeta(5) - 
9\,\zeta(4)\sp2}{\zeta (4)}\bigg\} \nsubseteq \Bbb Q, \quad k = 0,1.$$
\endproclaim
In the first part, I give a short survey of the C.\ S.\ Mejer functions,  
and then define and compute my auxiliary functions. In what follows, 
given $a \in \Bbb C$ and $M \subseteq \Bbb C, \; M \neq \emptyset,$ 
let 
$$
a + M  = M + a = \{w \in \Bbb C \, \colon \, w = a + z, \, z \in M \}.
$$ 

{\sl Acknowledgment.} I am obliged to Professor A.G.Aleksandrov and 
to Professor B.Z.Moroz for their help and support.

\heading{\sl \S\, 1. Short survey of the Mejer functions}\endheading 

In this work, I use the functions introduced and studied by C.\ S.\ Mejer 
in a long series of papers published during the decade 1936 -- 1946. The 
functions of Mejer can be defined as follows ([9], ch. 5). Let
$$
G \sb {p, q} \sp{(m, n)}\left(z\bigg\vert\matrix a\sb1,&\ldots,&a\sb p\\
a\sb1,&\ldots,&a\sb p\\
\endmatrix\right) = \frac{1}{2\pi i}\, \int \limits \sb L \frac {z\sp s \prod 
\limits \sb {j = 1} \sp m \Gamma(b\sb j - s)\prod \limits \sb {j = 1} 
\sp n \Gamma(1 - a\sb j + s)}{\prod \limits \sb {j = m + 1} \sp q 
\Gamma(1 -  b\sb j + s)\prod \limits \sb {j = n + 1} \sp p 
\Gamma(a\sb j - s)}\,ds, \tag{1.1} 
$$ 
where an empty product is, by definition, equal to 
$1, \; 0 \le  m \le  q, \; 0 \le  n \le  p;$  
the parameters $a\sb j \in \Bbb C, \; j = 1, \ldots, p,$ and 
$b \sb k \in \Bbb C, \; k = 1, \ldots, q,$ are chosen in such 
a way that none of the poles of the function $\Gamma(b \sb k - s), 
\; k = 1, \ldots, m,$ is equal to a pole of one of the functions 
$\Gamma(1 - a \sb j + s), \; j = 1, \ldots, n.$ There are 3 
possibilities to choose the curve $L.$ 

$(A)$ First, the curve $L = L\sb 0$ may be chosen to pass from $-i\infty$ 
to $+i\infty$ in such way that for any $k = 1, \ldots, m\,$ all the 
poles of the function $\Gamma(b\sb k - s)$ lie to the right of it 
and for any $j = 1, \ldots, n\,$ all the poles of the function 
$\Gamma(1 - a\sb j + s)$ lie to the left of it. The integral (1.1) 
is convergent in either of the following two cases: 

\smallskip
$(A1) \; \vert \arg(z) \vert < (m + n - \frac{p}{2} - \frac{q}{2})\pi ;$

\smallskip
$(A2) \; \vert \arg(z) \vert  \le  (m + n - \frac{p}{2} - \frac{q}{2})\pi$ 
and $\frac{p-q}{2} + \Re\Delta\sp* < -1 ,$ where 
$$\Delta \sp * = \sum \limits \sb {k = 1} \sp q b\sb k  - 
\sum \limits \sb {k = 1} \sp q a\sb j. \tag{1.2}$$

$(B)$ Second, the curve $L = L\sb 1$ may be chosen to pass from 
$+i\infty$ to $+i\infty,$ encircling each of the poles of the functions 
$\Gamma(b\sb k - s), \; k = 1, \ldots, m,$ in the negative direction, 
but not including any of the poles of the functions 
$\Gamma(1 - a \sb j + s), \, j = 1, \ldots, n.$ The integral (1.1) 
is convergent in each of the following three cases: 

\smallskip
$(B1) \;  p < q ;$

\smallskip
$(B2) \; 1 \le p \le  q$ and $\vert z \vert  < 1 ;$ 

\smallskip 
$(B3) \; 1 \le p \le  q, \; \vert z \vert  \le  1$ and 
$\Re\Delta\sp* < -1 .$ 

\medskip
$(C)$ Third, the curve $L = L\sb 2$ may be chosen to pass from 
$-\infty$ to $-\infty,$ encircling each of the poles of the functions 
$\Gamma(1 - a \sb j + s), \; j = 1, \ldots, n,$ in the positive 
direction, but not including any of the poles of the functions 
$\Gamma(b\sb k - s), \; k = 1, \ldots, m.$ 
The integral (1.1) is convergent in each of the following 3 cases: 

\smallskip
$(C1) \; q < p;$

\smallskip
$(C2) \; 1 \le q \le  p$ and $\vert z \vert  > 1;$ 

\smallskip
$(C3) \; 1 \le q \le  p, \; \vert z \vert  \ge  1$ and 
$\Re\Delta\sp* < -1.$

\smallskip
If both conditions (A) and (B) (respectively (A) and (C)) are satisfied, 
then the result does not depend on whether the curve $L$ is defined as in 
(A) or as in (B) (respectively in (C)). 

Let $g$ denote the integrand of the integral (1.1), let $G$ denote 
the integral (1.1) with $L = L\sb k,$ where $k=1,2,$ and let $S\sb k$ be  
the set of all the unremovable singularities of $g$ encircled by $L\sb k.$ 
If one of the conditions (B1) -- (B3) (respectively (C1) -- (C3)) holds 
for $k = 1$ (respectively for $k = 2),$ then 
$$G = (-1)\sp k \sum \limits \sb {s \in S\sb k}\Res(g;s), \tag{1.3}$$ 
where $\Res(g;s)$ stands for the residue of the function $g$ at the point 
$s.$ Let us prove these assertions. 
Let $\sigma \sb 0 > 0, \; \tau \sb 0 > 0,$ and let 
$$\min(\sigma \sb0, \tau \sb0) \ge 1 + \sum \limits \sb {j = 0} 
\sp p 2(\vert a\sb j\vert  + 1)  + \sum \limits \sb {k = 0} 
\sp q 2(\vert b\sb k\vert  + 1). \tag{1.4}$$

The curve $L\sb 0$ may be chosen in such a way that, except for
a compact piece, it coincides with a part of the imaginary axis 
$\vert \Im{s} \vert \ge \tau \sb 0;$ the curves $L\sb 1$ and $L\sb 2 $ may 
be chosen so that, except for a compact part, the curve $L\sb 1$ 
coincides with the union of the two rays $\sigma  \pm  i\tau \sb0, \; 
\sigma \ge \sigma \sb 0,$ and the curve $L\sb 2,$ except for a compact part, 
coincides with the union of the two rays $\sigma  \pm  i\tau \sb0, \; 
\sigma  \le  -\sigma \sb 0.$ 

Let $r\sb 0 \ge \sigma \sb0  + \tau \sb 0$ be chosen big enough, so that 
the specified compact parts of $L\sb0, \, L\sb 1$ and $L\sb 2$ lie inside 
the disk $K\sb 0 = \lbrace s \in \Bbb C \,\colon\, \vert s \vert \le 
r\sb 0\rbrace.$ 

Let $u \ge  0, \; \eta (u) = \int \limits \sb 0 \sp u 
(\lbrack t\rbrack - t + \frac12)dt;$ then $0 \le  \eta (u) \le  1/8.$ 
Let  
$$
\overline{D} \sb0  = \lbrace s \in 
\Bbb C \,\colon\,  \vert s \vert \ge r\sb 0\rbrace \backslash  \lbrace 
s \in \Bbb C \,\colon\, \Re{s} < 0, \; \vert \Im{s} \vert  < \tau \sb 0 
\rbrace.
$$
 
If $z = x + iy, \; \vert y \vert \ge \tau \sb 0, $ then the integral 
in the complex Stirling formula ([10], ch. 4) 
$$ 
\log \Gamma (z) = (z - \frac12)\log{z} - z + \frac12 \log(2\pi) + 
\int \limits \sb 0 \sp \infty \frac {\eta (u)}{(u + z)\sp 2}du 
$$
can be estimated as follows: 
$$
\bigg\vert \int \limits \sb 0 \sp \infty \frac {\eta (u)}{(u + z)\sp 2}du 
\bigg\vert  \le   \int \limits \sb 0 \sp \infty  \bigg\vert  
\frac {\eta (u)}{(u + z)\sp 2} \bigg\vert du =
$$
$$
\int \limits \sb 0 \sp \infty  \bigg\vert\frac {\eta (u)}{(u + x)\sp 2 
+ y\sp 2}\bigg\vert du < \int \limits \sb {- \infty} \sp {+\infty} 
\bigg\vert \frac {\eta (u)}{(u + x) \sp 2 + y\sp 2}\bigg\vert du 
< \frac{\pi}{8\tau \sb 0} = O(1).
$$ 

If $z = x + iy$ and $x \ge \sigma \sb 0, $ then 
$$
\bigg\vert \int \limits \sb 0 \sp \infty \frac {\eta (u)}{(u + z)\sp 2}du 
\bigg\vert  \le   \int \limits \sb 0 \sp \infty  \bigg\vert  
\frac {\eta (u)}{(u + z)\sp 2} \bigg\vert du =
$$
$$
\int \limits \sb 0 \sp \infty \bigg\vert\frac {\eta (u)}
{(u + x)\sp 2 + y\sp 2}\bigg\vert du \le  
\int \limits \sb 0 \sp \infty \bigg\vert \frac {\eta (u)}
{(u + x)\sp 2} \bigg\vert du < \frac{1}{8\sigma \sb 0} = O(1).
$$ 
Therefore 
$$
\log \Gamma (s + c) = (s + c - \frac12)(\log s - 1) + O(1) \tag{1.5}
$$
in $\overline{D}\sb 0,$ for any constant $c$ with $2\vert c \vert \le 
\min\{\sigma \sb0, \tau \sb 0\}.$ 

Let $s = \sigma  + i\tau  = \vert s\vert \,e \sp {i\psi}$ with 
$\vert \psi \vert  \le \pi, \; \phi  = \arg z,$ then 
$$ 
\Re\log\Gamma(s + c) = \big(\sigma + \Re(c) - \frac12\big)(\log\vert
s\vert - 1) - \tau \psi  + O(1). \tag{1.6}
$$ 
One may write $g(s) = g\sb 1(s)g\sb 2(s)$ with 
$$
g\sb 1(s) =  z\sp s \prod \limits \sb {j = 1} \sp p 
\Gamma(1 -a\sb j + s)\bigg/\prod \limits \sb {j = 1} \sp q 
\Gamma (1 - b\sb j + s), \tag{1.7}
$$
$$
g\sb 2 (s) = \prod \limits \sb{j = n + 1}\sp p \frac{\sin(\pi(a\sb j - s))}\pi 
\prod \limits \sb {j = 1} \sp m \frac \pi {\sin(\pi (b\sb j -s))} . \tag{1.8}
$$ 
Since $\tau \psi \ge 0,$ it follows from (1.6) that 
$$
g\sb 1(s) = O(1)\cdot \vert z\vert\sp \sigma e\sp {\vert \phi \vert \vert \tau \vert}(\vert s \vert/e) \sp 
{((p - q)(\sigma  + \frac12) + \Re\Delta\sp*)} e\sp {-\tau (p - q)\psi} = 
\tag{1.9}
$$
$$
= O(1)\cdot \vert z\vert\sp \sigma(\vert s \vert/e) \sp {((p - q)(\sigma  + \frac12) 
+ \Re\Delta\sp*)} e\sp {\vert \tau \vert (p - q)(\frac{\pi}{2} - \vert 
\psi \vert)}e\sp {\vert \tau \vert(\phi - (p - q)\frac{\pi}{2})} 
$$
in $\overline{D}\sb 0.$ 
In view of our choice of $\tau \sb 0,$ it follows that 
$$
\vert \sin(\pi(a\sb j - s)) 
\vert  = O(1)\cdot e\sp {\pi \vert \tau \vert}, \quad \vert 
\sin(\pi(b\sb k - s)) \vert \sp {-1} = 
O(1)\cdot e\sp {-\pi \vert \tau \vert}
$$
for $s \in \overline{D}\sb 0, \; j = 1, \ldots, p$ and $k  = 1, \ldots, q,$ 
and consequently 
$$ 
g\sb 2(s) = O(1)\cdot e\sp {\pi (p - m - n)\vert \tau \vert}\,. \tag{1.10}
$$
Therefore, according to (1.9) and (1.10) the following equality 
$$
\vert g(s) \vert  = O(1)\cdot \vert z \vert \sp \sigma e\sp {\vert \phi
\vert \vert \tau \vert} (\vert s \vert/e) \sp {((p - q)(\sigma  + \frac12)
+ \Re\Delta\sp*)} e\sp {-\vert \tau \vert (p - q)\vert \psi \vert}e\sp
{\pi (p -  m - n)\vert \tau \vert} = \tag{1.11}
$$
$$
= O(1)\cdot \vert z \vert \sp \sigma (\vert s \vert/e) \sp {((p - q)
(\sigma  + \frac12) + \Re\Delta\sp*)}\cdot
e\sp {((p - q)(\frac{\pi}{2} - \vert \psi \vert) + \vert \phi \vert + \pi 
(\frac{p+q}{2} - m- n))\vert \tau \vert}
$$ 
holds in $\overline{D}\sb 0.$  
According to (1.11), on the part of the curve $L\sb 0$  
lying on the rays $\sigma  = 0, \, \vert \tau \vert  \ge  \tau \sb 0$ 
we have the equality 
$$
\vert g(s)\vert  = O(1)\cdot(\vert \tau \vert/e) \sp {(\frac{p-q}{2} + 
\Re\Delta\sp*)} 
e\sp {\vert \tau \vert(\vert \phi \vert + \pi(\frac{p+q}{2} -  m - n))};
$$ 
if the condition (A1) holds, then 
$$
\vert g(s)\vert  =\vert \tau \vert \sp {O(1)}e\sp {\vert \tau 
\vert(\vert \phi \vert + \pi(\frac{p+q}{2} -  m - n))},
$$ 
where $\vert \tau \vert \ge \tau \sb 0,$ and the integral (1.1) 
is convergent; if the condition (A2) holds, then 
$$
\vert g(s)\vert  =  O(1)\cdot (\vert \tau \vert/e) \sp {(\frac{p-q}{2} + 
\Re\Delta\sp*)},
$$
where $\vert \tau \vert \ge \tau \sb 0,$ and the integral (1.1) 
is convergent. Let the number $\sigma \sb 0$ be chosen so that 
$$
\sigma \sb 0 > 3 + \bigg\vert \frac{p-q}{2} + \Re\Delta\sp*\bigg\vert. \tag{1.12}
$$ 
If one of the condition (B), (C) holds, then 
$$
(p - q)\sigma = -\vert(p - q)\vert\vert\sigma\vert, \tag{1.13}
$$
$$
(p - q)(\frac{\pi}{2} -\vert \psi \vert) = -\vert p - q\vert \vert 
\frac{\pi}{2} -\vert \psi \vert \vert, \tag{1.14}
$$
$$
\sigma(\log\vert z\vert) = -\vert \sigma \vert \vert 
\log\vert z\vert\vert \tag{1.15}
$$ 
in $\overline{D}\sb 0,$ and (1.9), (1.11) takes respectively the form 
$$
\vert g\sb1(s)\vert  = O(1)\cdot e\sp {-\vert \sigma\vert  
\vert \log \vert z\vert\vert}(\vert s\vert/e)\sp{-\vert \sigma \vert  
\vert p - q\vert + \frac{p-q}{2} + \Re\Delta\sp*}\times \tag{1.16}
$$
$$
e\sp {-\vert p - q\vert  \vert\frac{\pi}{2} - \vert \psi \vert  
\vert \vert \tau \vert} e \sp  {(\vert \phi \vert  - 
(p - q) \frac{\pi}{2} )\vert \tau \vert }
$$  
$$
\vert g(s)\vert  = O(1)\cdot e\sp {-\vert \sigma\vert  \vert 
\log \vert z\vert\vert}(\vert s\vert/e)\sp{-\vert \sigma \vert  \vert 
p - q\vert + \frac{p-q}{2} + \Re\Delta\sp*}\times \tag{1.17}
$$
$$
e\sp {(-\vert p - q\vert\vert  \frac{\pi}{2} - \vert \psi \vert \vert  
+ \vert \phi \vert + \pi (\frac{p+q}{2} - m- n))\vert \tau \vert}.
$$ 
If one of conditions (B1), (C1) holds and $\vert\sigma\vert \ge 
\sigma \sb 0,$ then $\vert p - q\vert \ge 1$ and, in view of 
(1.12) and (1.17), we have 
$$
\vert g\sb1(s)\vert  = O(1)\cdot e\sp {-\vert \sigma\vert  \vert 
\log\vert z\vert\vert}(\vert s\vert/e)\sp{-3}e\sp {-\vert p - q\vert  
\vert\frac{\pi}{2} - \vert \psi \vert  \vert}e \sp  {(\vert \phi \vert  - 
(p - q) \frac{\pi}{2} )\vert \tau \vert}, \tag{1.18}
$$ 
$$
\vert g(s)\vert  = O(1)\cdot e\sp {-\vert \sigma\vert  \vert 
\log\vert z\vert\vert}(\vert s\vert/e)\sp {-3} 
e\sp {(-\vert p - q\vert\vert \frac{\pi}{2} - \vert \psi \vert  \vert  + 
\vert \phi \vert + \pi (\frac{p+q}{2} - m- n))\vert \tau \vert}  \tag{1.19}
$$
in $\overline{D}\sb 0.$ Therefore if one of conditions (B1), (C1) holds, 
then on the rays $\tau = \tau \sb 0,$ $\vert \sigma \vert  \ge 
\sigma \sb 0$ we have the equality 
$$
\vert g(s)\vert  = O(1)\cdot(\vert \sigma \vert/e)\sp {-3}, \tag{1.20}
$$ 
and the integral (1.1) is convergent; if one of conditions 
(B2), (C2) holds, then on the rays $\tau = \tau \sb 0, \; 
\vert \sigma \vert  \ge \sigma \sb 0$ with $\vert 
\log \vert z\vert\vert > 0$ we have the equality 
$$
\vert g(s)\vert  = \vert \sigma \vert \sp {O(1)}\cdot 
e\sp {-\vert \sigma\vert \vert \log \vert z\vert\vert}, \tag{1.21}
$$ 
and the integral (1.1) is convergent. Finally, if one of  
condition (B3), (C3) holds, then, in view of (1.17), 
$$
\vert g(s)\vert  = O(1)\cdot \vert\sigma\vert\sp{\Re\Delta\sp*} 
\tag{1.22}
$$ 
on the rays $\tau = \tau \sb 0, \; \vert \sigma \vert  > \sigma \sb 0,$ 
and the integral (1.1) is convergent. For any $r \ge r\sb 0$ and 
$k = 0, 1, 2$ let us denote by $L\sb {r, k}$ the part of the curve 
$L\sb k$ contained in the disk $K\sb r = \{ s \in  \Bbb C \,\colon\, 
\vert s \vert \le r \}$ and directed as the curve $L\sb k.$ Let 
further $C\sb {r, k, j},$ where $k = 1, 2, \, j = 0, 1,$ be the smallest 
counter-clockwise directed arc of the circle $\vert s \vert  = r$ 
connecting the point $ir(-1)\sp j$ with the point 
$i\tau \sb 0 (-1)\sp j - (-1)\sp k\sqrt{r\sp 2 - \tau \sb 0\sp2}.$ 
Finally, let $L\sb {r,k} \sp *,$ where $k = 1,2,$ be the 
positively directed curve consisting of the curves 
$L\sb{0, r}, \, L\sb {k, r}, \, C\sb {r, k, 0}$ and $C\sb {r, k, 1}.$ 
Clearly, 
$$
\int \limits \sb {L\sb k} g(s) = \lim \limits \sb 
{r\rightarrow \infty}\int \limits \sb {L\sb {r,k}} g(s)ds \tag{1.23}
$$
for $k = 0, 1, 2.$ 
The function $g(s)$ has no singular points in the domain bounded 
by the curve $L\sb k \sp *.$ Therefore 
$$
\int \limits \sb {L\sb {r, k}}\hskip-4pt g(s)ds - (-1)\sp k 
\hskip-6pt \int \limits \sb {C\sb {r, k, 0}}\hskip-6pt g(s)ds - 
\int \limits \sb {L\sb {r, 0}}\hskip-4pt g(s)ds - 
(-1)\sp k\hskip-6pt\int \limits \sb {C\sb {r, k, 1}}\hskip-6pt g(s)ds 
\; = \int \limits \sb {L\sb k \sp *} g(s)ds = 0\,. \tag{1.24}
$$ 
For $k = 1, 2$ and $j = 0, 1,$ let $\check C\sb {r, k, j}$ be 
the counter-clockwise directed part of the arc $C\sb {r, k, j}$ 
lying in the strip $\vert \sigma \vert \le \sigma \sb 0,$ and let 
$\hat C\sb {r, k, j}$ be the counter-clockwise directed part of 
the arc $C\sb {r, k, j}$ lying outside this strip. 
If $p \ne q$ and one of conditions (B1), (C1) holds, 
then according to (1.19)
$$
\vert g(s)\vert  = O(1)\cdot (r/e)\sp {-3} \tag{1.25}
$$
on the $\hat C\sb {r, k, j}$ for $k = 1, 2$ and $j = 0, 1,$ 
and, because the length of the arc $\hat C\sb {r, k, j}$ is not 
bigger than $\pi r/2,$ it follows that 
$$
\lim \limits \sb {r \rightarrow  \infty}  \int \limits 
\sb {\hat C\sb {r, k, j}} g(s)ds = 0. \tag{1.26}
$$ 
Let $\gamma = 4\sp{-1}\min\{\vert \vert \phi \vert + 
\pi (p\ +\ q)/2 - m - n\vert,\vert \log\vert z \vert \vert\}.$
In view of (1.17), if $p=q$ and one of conditions (B), (C) 
holds, then 
$$
\vert g(s)\vert  = O(1)\cdot e\sp {-\gamma r} (r/e)\sp {\Re\Delta\sp*} 
\tag{1.27}
$$ 
on the $\hat C\sb {r, k, j}$ for $k = 1, 2$ and $j = 0, 1;$
now, if the condition (A) is satisfied, then the equality $\gamma = 0$
implies the inequality $\Re(\Delta\sp *)<-1,$ and because 
the length of the arc $\hat C\sb {r, k, j}$  is not bigger than 
$\pi r/2,$ relation (1.26) still holds. Moreover, 
$$
\vert \tau \vert  = \sqrt{r\sp 2 - \sigma \sb 0\sp 2}  = 
r - O(1), \quad \vert \psi \vert  = \frac{\pi}{2} - 
O(\arcsin(\sigma\sb 0 / r)) = \frac{\pi}{2} - O(1/r)
$$ 
on the $\check C\sb {r, k, j}$ for $k = 1, 2, \, j = 0, 1.$ 
Therefore if one of conditions (B), (C) holds, then according to 
(1.17)
$$
\vert g(s)\vert  = O(1)\cdot (r/e)\sp{\frac{p-q}{2} + \Re\Delta\sp*}
e\sp { (\vert \phi \vert + \pi (\frac{p+q}{2} - m- n))\vert \tau \vert} = 
\tag{1.28}
$$
$$
O(1)\cdot (r/e)\sp{ \frac{p-q}{2} + \Re\Delta\sp*}e\sp {(\vert \phi \vert 
+ \pi (\frac{p+q}{2} - m- n))r}
$$ 
on the curve $\check C\sb {r, k, j}.$ 
Consequently, if the condition (A1) and one of the condition (B), (C) 
hold, then
$$
\lim \limits \sb {r \rightarrow  \infty}  \int \limits \sb 
{\check C\sb {r, k, j}} \hskip-4pt g(s)ds = 
\lim \limits \sb {r \rightarrow  \infty} O(1)\cdot e\sp{(\vert \phi 
\vert + \pi (\frac{p+q}{2} - m- n))r}r\sp{O(1)} = 0 \tag{1.29} 
$$ 
because in this case $\vert\phi\vert + \pi(\frac{p+q}{2} - m- n) < 0.$ 
If the condition (A2) holds, then 
$$
\lim \limits \sb {r \rightarrow  \infty}  
\int \limits \sb {\check C\sb {r, k, j}} \hskip-4pt g(s)ds = 
\lim \limits \sb {r \rightarrow  \infty}O(1)\cdot r\sp{ \frac{p-q}{2} + 
\Re\Delta\sp* + 1} = 0. \tag{1.30}
$$ 
We have thus proved that if both of the convergence conditions 
(A) and (B) hold, then 
$$
\int \limits \sb {\check L\sb 0} g(s)ds = 
\int \limits \sb {L\sb 1} g(s)ds;
$$ 
and if both of the convergence conditions (A) and (C) hold, then 
$$
\int \limits \sb {\check L\sb 0} g(s)ds = 
\int \limits \sb {L\sb 2} g(s)ds.
$$ 
Now I shall prove the equality (1.3). Let us choose $\sigma\sb 1$ 
so that it is bigger than $r\sb0$ and lies outside of the countable set
$$
\bigcup \limits \sb {\mu  = 0} \sp 1 \bigcup \limits \sb {j = 1} 
\sp p \, \big((-1)\sp \mu \Re(a\sb j) + \Bbb Z\big) \, 
\bigcup \bigcup \limits \sb {\nu  = 0} \sp 1 \, 
\bigcup \limits \sb {k = 1} \sp q \, \big((-1)\sp \nu 
\Re(b\sb k) + \Bbb Z\big).
$$
For $k = 1, 2$ and  $n\in \Bbb N - 1 ,$ 
let $r(n) = \sqrt{\tau \sb 0\sp 2 + (\sigma \sb 1 + n)\sp 2},$  
let $c\sb {n, k}$ denote the segment 
$$
[-(-1)\sp k(\sigma\sb1 + n) - i\tau \sb 0, 
-(-1)\sp k(\sigma\sb1 +  n) + i\tau \sb 0],
$$ 
and let $\widetilde L\sb {n, k}$ be the counter-clockwise directed 
curve consisting of the curve $L\sb {r(n), k}$ 
and the segment $c\sb {n, k};$ finally, let $\widetilde S\sb {n, k}$ 
denote the set of all the unremovable singularities encircled by the 
curve $\widetilde L\sb {n, k}.$ It is clear that 
$$
\frac{1}{2\pi i} \int \limits \sb {\widetilde S\sb {n, k}} g(s)ds = 
\sum \limits \sb {s  \in \widetilde S\sb {n, k}} \Res(g ; s),
$$
where $k=1,2.$ The equality (1.3) will be proved if we establish that 
for each $k = 1, 2$ the equality
$$
\lim \limits \sb{n \rightarrow  \infty}   \int \limits \sb 
{c\sb{n, k}} g(s)ds = 0 \tag{1.31}
$$ 
holds. First, consider the case $k = 1.$ In view of the choice of 
$\sigma \sb 1,$ the function $g\sb 2(s)$ is continuous on the 
segment $c\sb {0, k},$ and therefore there is a constant 
$M\sb 0$ with $\vert g\sb 2(s)\vert \le  M\sb 0$ for $s\in c\sb {0,k},$ 
where $k=1,2.$ Since the period of the function $\vert g\sb 2 (s)\vert$ 
is equal to one, we have 
$$
\vert g\sb 2(s)\vert  \le  M\sb 0 \tag{1.32}
$$ 
for $s \in c\sb {n, k}, \; k = 1, 2,$ and $n \in  \Bbb N.$ 
Let us first consider the case $k = 1.$ 
If condition (B) holds, then (1.31) follows from (1.16). 
Suppose now that condition (C) is satisfied. Write 
$$
g\sb 1 (s) = g\sb 3(s)g\sb 4(s) \tag{1.33}
$$
with 
$$
g\sb 3(s) =  z\sp s \prod \limits \sb {j = 1} \sp q 
\Gamma (b\sb j - s)\bigg/\prod \limits \sb {j = 1} \sp p 
\Gamma(a\sb j - s), \tag{1.34}
$$
$$
g\sb 4 (s) = \bigg(\prod \limits \sb {j = 1} \sp p 
\frac {\sin(\pi (a\sb j - s))}{\pi} \prod \limits \sb {j = 1} \sp q 
\frac{\pi}{\sin(\pi (b\sb j -s))}\bigg)\sp{\!-1} . \tag{1.35}
$$ 
Since the period of the function $\vert g\sb 4 (s)\vert$ is equal 
to one, there is a constant $M\sb 2$ such that 
$$
\vert g\sb 4(s)\vert  \le  M\sb 2 \tag{1.36}
$$ 
for $s \in c\sb {n, 2}, \; n \in  \Bbb N - 1 ,$ in view of the choice 
of $\sigma \sb 1.$ Because $\vert \tau \vert \le \tau \sb 0$ 
if $s \in c\sb {n, 2}$ and $-c\sb {n, 2}$ lies in $\overline{D}\sb 0, \; 
n \in \Bbb N - 1 ,$ it follows from (1.6) that 
$$
g\sb 1(s) = O(1)\!\cdot\! \vert z\vert\sp \sigma 
\big(\vert s \vert /e \big) \sp {\sum \limits \sb {j = 1} \sp q 
\Re(b\sb j - \sigma - \frac12) - \sum \limits \sb {j = 1} \sp p 
\Re(b\sb j - \sigma - \frac12)} = 
$$
$$
= O(1)\!\cdot\! \vert z\vert\sp \sigma \big(\vert s \vert /e\big) 
\sp{((p - q)(\sigma  + \frac12) + \Re\Delta\sp*)} 
$$ 
with $\vert s \vert \ge -\sigma = \sigma\sb1 + n.$ 
Therefore the equality (1.31) holds also for $k = 2.$ 
This proves (1.3). 

\heading {\sl \S\, 2. My auxiliary functions}\endheading

I shall work with the set $\Omega \sp *$ consisting of all the points 
$z  \in \Bbb C$ for which 
$$
\vert z \vert  \ge  1 \quad\text{and}\quad 
-\frac{3 \pi}{2} < \arg(z) \le  \frac{\pi}{2} \,. \tag{2.1}
$$
Thus $\log(-z) = \log{z} - i\pi$ for $\Re(z) > 0.$ 
Let $ \Delta  \in  \Bbb N , \; 1 < \Delta, \; 
\delta \sb 0  = 1 / \Delta,$ 
$$
\gamma \sb 1 = (1 - \delta \sb 0 ) / (1 + \delta \sb 0 ), \quad 
d \sb {l} = \Delta +(-1)\sp l, \quad l = 1, 2.
$$ 
To introduce the first of my auxiliary function $f\sb 1 (z, \nu),$ I use 
the auxiliary set 
$$
\Omega \sp {(h)} = \lbrace z \in \Bbb C \,\colon\, \vert z\vert \le 1, \; 
-\frac{3\pi}{2} < \phi = \arg(z) \le \frac{\pi}{2}\rbrace. 
$$ 
I shall prove that, for each $\nu \in \Bbb N,$ the function 
$f\sb 1 (z, \nu)$ belongs to $\Bbb Q[z];$ therefore using 
the principle of analytic continuation we may regard it as being 
defined in $\Bbb C,$ and consequently, in $\Omega \sp * .$ 
For $\nu \in \Bbb N,$ let 
$$
f\sb 1 (z, \nu) =  \tag{2.2}
$$
$$
= (-1) \sp {\nu (\Delta + 1)}G \sb {6, 6} \sp{(1, 3)}
\left(z\bigg\vert\matrix -\nu d\sb1, &\!\! -\nu d\sb1, &\!\! -\nu d\sb1, 
&\!\! 1+\nu d\sb2, &\!\! 1+\nu d\sb2, &\!\! 1+\nu d\sb2\\  
0, &\!\! 0, &\!\! 0, &\!\! \nu, &\!\! \nu, &\!\! \nu\\
\endmatrix\right) 
$$
$$
= -(-1) \sp {\nu (\Delta + 1)}\frac{1}{2\pi i} \int \limits \sb 
{L \sb 1} g \sb {6, 6} \sp{(1, 3)}(s)ds,
$$ 
where 
$$
g \sb {6,6} \sp{(1,3)} (s) = z\sp s\Gamma (-s)\Gamma(1 + s) 
\sp{-2}\big(\Gamma(1 + \nu d\sb 1 + s)/(\Gamma (1 - \nu +s)\Gamma 
(1 + \nu d \sb 2 - s))\big)\sp3,
$$  
and the curve $L\sb1$ passes from $+\infty$ to $+\infty$ encircling 
the set $\Bbb N - 1 $ in the negative direction, but not including any 
point of the set $- \Bbb N.$  Here $p = q = 6, \, m = 1, \, 
n = 3,$ $a\sb 1 = a\sb 2 = a\sb 3 = -\nu d\sb 1, \, 
a\sb 4 = a\sb 5 = a\sb 6 = 1 + \nu d \sb 2, \, 
b\sb 1 = b\sb 2 = b \sb 3 = 0,$ $b\sb 4 = b\sb 5 =  b\sb 6 = \nu,$  
$\Delta \sp * = -3\nu - 3$
and, since we take $\vert z \vert \le  1,$ convergence conditions (B2) 
and (B3) hold. To compute the function $f\sb 1 (z, \nu),$ we use 
formula (1.3) and the well-known formula 
$$
\Gamma(s) = \Gamma(s + l) \prod \limits \sb {k = 1} 
\sp l (s + l - k)\sp {-1} \tag{2.3}
$$ 
with $l \in \Bbb N.$ The set of unremovable singular points 
of the function $g \sb {6, 6} \sp{(1, 3)} (s),$ which are encircled by 
the curve $L\sb 1,$ consists of the points 
$s = \nu, \ldots, \nu d\sb2,$ all these points are poles of the first 
order, and, for $k = 0, \ldots, \nu \Delta,$ the following equality holds: 
$$
\Res(g \sb {6, 6} \sp{(1, 3)};\nu  + k) = 
-(-z)\sp {\nu  + k}((\nu  + k)!)\sp {-3}((\nu \Delta + k)!)\sp 3(k!)\sp{-3} 
((\nu \Delta - k)!)\sp{-3} =
$$ 
$$
= (-z)\sp{\nu  + k}\big((\nu d\sb 1)!/(\nu \Delta)!\big)\sp 3 
\binom {\nu \Delta} k \sp{\! 3}\binom{\nu\Delta  + k}{\nu d\sb 1}\sp{\!3}.
$$ 
The function $f\sb 1 (z, \nu)$ is equal to a finite sum 
$$
f\sb 1(z, \nu) = \big((\nu d\sb 1)!/(\nu \Delta)!\big)\sp 3 
z\sp \nu(-1)\sp {\nu\Delta}\sum \limits \sb {k = 0} \sp{\nu\Delta}(-z) 
\sp k \binom {\nu \Delta  }k\sp{\!3} \binom {\nu \Delta + k  }{\nu d \sb 1} 
\sp{\!3}. \tag{2.4}
$$ 
Therefore, as it has been already remarked, using the principle of 
analytic continuation we may regard it as being defined in 
$\Bbb C$ and consequently in $\Omega \sp * .$ 

Let $\Omega \sp * \sb0  = \lbrace z\in \Omega \sp * 
\,\colon\, \Re(z) > 0\rbrace.$ 
If $z \in  \Omega \sp * \sb 0,$ then 
$$-\frac{3\pi}{2} < \phi = \arg(-z) = \arg(z) - \pi < -\frac{\pi}{2},$$ 
and therefore $-z \in  \Omega \sp *.$ 
Now, let me introduce my second auxiliary function defined 
for $z \in  \Omega \sp * \sb 0.$  

Let
$$
f\sb 2(z, \nu) = -(-1) \sp {\nu \Delta}\,\frac{1}{2\pi i}
\int \limits \sb{L \sb 2} g \sb {6, 6} \sp{(4, 3)}(s)ds =  
$$
$$ 
= -(-1) \sp {\nu \Delta} G \sb {6, 6} \sp{(4, 3)}
\left(- z \bigg\vert\matrix -\nu d\sb1, &\!\! -\nu d\sb1, 
&\!\! -\nu d\sb1, &\!\! 1+\nu d\sb2, &\!\! 1+\nu d\sb2, &\!\! 1+\nu d\sb2\\
0, &\!\! 0, &\!\! 0, &\!\! \nu, &\!\! nu, &\!\! \nu\\
\endmatrix\right),
$$ 
where $z \in  \Omega \sp * \sb0, \, \nu \in \Bbb N,$ and 
$$g \sb {6, 6} 
\sp {(4, 3)} = g \sb {6, 6} \sp{4, 3)} (s) = 
(-z)\sp s \Gamma(-s)\sp{3}\Gamma (-s + \nu) 
\Gamma(1 - \nu + s)\sp{-2}\Gamma (1 + \nu d\sb 1 + s)\sp3
\Gamma (1 + \nu d \sb 2 - s)\sp{-3}{\!\!,}
$$ 
and the curve $L\sb 2$ passes from $- \infty$ to $- \infty,$ 
encircling the set $-\Bbb N$ in the positive direction but no 
point in the set $\Bbb N - 1 .$ Here $p = q = 6, \, m = 4, \, n = 3,$ 
$a\sb 1 = a\sb 2 = a\sb 3 = -\nu d\sb 1, \, 
a\sb 4 = a\sb 5 = a\sb 6 = 1 + \nu d \sb 2, \, 
b\sb 1 = b\sb 2 = b \sb 3 = 0,$
$b\sb 4 = b\sb 5 =  b\sb 6 = \nu, \, \Delta \sp * = -3\nu - 3;$ 
since now $\vert -z \vert \ge  1,$ convergence conditions (C2) and (C3) 
are satisfied. To compute the function $f\sb 2 (z,\nu),$ we use formula (1.3). 
The set of all the unremovable singular points of the function 
$g \sb {6, 6} \sp{(4, 3)} (s),$ encircled by the curve $L\sb 2,$ 
consists of the points $s  = -1 - \nu d\sb 1 - k$ with $k \in  \Bbb N - 1 ;$ 
each of these points is a pole of the first order. Therefore making use 
of (2.3) one obtains 
$$ 
\Res(g \sb {6, 6} \sp{(4, 3)}; -1 - \nu d\sb 1 - k) =  
$$
$$
= (-z)\sp { -1 - \nu d\sb 1 - k} ((\nu d\sb 1 + k)!)
\sp 3((\nu \Delta + k)!)\sp3(-1)\sp k  
(k!)\sp{-3}((1 + 2\nu \Delta + k)!)\sp{-3} =
$$ 
$$
= (-1)\sp {1 + \nu  d \sb 1}z\sp {-(1 + \nu  d \sb 1 + k)}
\left(\frac {\prod \limits \sb {j = 1} \sp {\nu \Delta  - \nu}(1 + \nu 
\Delta  - \nu + k - j)}{\prod \limits \sb {j = 0} \sp {\nu \Delta} 
(1 + \nu \Delta  + k + j)}\right)\sp{\! 3}\!,
$$ 
$$
f\sb 2(z, \nu) = \sum \limits \sb {k = 0} \sp {+ \infty} z
\sp {-(1 + \nu  d \sb 1 + k) }\left(\frac {\prod \limits \sb {j = 
1} \sp {\nu \Delta  - \nu}(1 + \nu \Delta  - \nu  + k - j)}
{\prod \limits \sb {j = 0} \sp {\nu \Delta}(1 + \nu \Delta  + 
k + j)}\right)\sp{\! 3}\!. \tag{2.5}
$$ 
Let $a \in  \Bbb N - 1 , \, b  \in \Bbb N + a,$ and 
$$
R(a;b;t) = \frac {b!}{(b - a)!}\prod \limits \sb {\kappa = 
a+1 } \sp b \! (t - \kappa) \prod \limits \sb {\kappa = 0 } \sp b 
\frac {1} {(t + \kappa)}, \qquad R \sb {0} (t;\nu) = R(\nu;\nu\Delta ;t). 
\tag{2.6}
$$ 
Let $t = 1 + \nu \Delta + k$ with $k \in  \Bbb N - 1 ;$ in view 
of (2.5), it follows that 
$$
f\sb 2 (z, \nu )\big((\nu \Delta )!/(\nu d\sb1)!\big)\sp 3 = 
(-1) \sp \nu \!\! \sum \limits \sb {t=\nu \Delta + 1} \sp \infty 
R \sb 0(t;\nu)\sp3 z\sp{-t + \nu} . \tag{2.7}
$$ 
Since $R \sb {0} (t;\nu) = 0$ for $t = \nu + 1, \ldots, \nu \Delta,$ 
we have 
$$
f\sb 2 (z, \nu )\big((\nu \Delta )!/(\nu d\sb1)!\big)\sp 3 =  
(-1) \sp \nu \!\!\sum \limits \sb {t=\nu + 1} \sp \infty R \sb 0(t;\nu)
\sp3  z\sp{-t + \nu} . \tag{2.8}
$$ 
Let 
$$
f\sb 3 = (-1) \sp {\nu (\Delta + 1)}\frac{1}{2\pi i} 
\int \limits \sb{L\sb2} g \sb{6, 6} \sp{(5, 3)}(s)ds =
$$ 
$$
(-1) \sp {\nu (\Delta + 1)}G \sb {6, 6} \sp{(5, 3)}\left(z \bigg\vert 
\matrix -\nu d\sb1, &\!\! -\nu d\sb1, &\!\! -\nu d\sb1, &\!\! 1+\nu d\sb2, 
&\!\! 1+\nu d\sb2, &\!\! 1+\nu d \sb2\\
0, &\!\! 0, &\!\! 0, &\!\! \nu, &\!\! \nu, &\!\! \nu\\
\endmatrix\right),
$$
where $z \in  \Omega \sp * \sb0, \,  \nu \in \Bbb N,$ and 
$$ 
g \sb {6, 6} \sp{(5, 3)} = g \sb {6, 6} \sp{(5, 3)} (s) = 
z\sp s \Gamma (-s)\sp{3}\Gamma (\nu - s)\sp 2 
\Gamma (1 - \nu + s)\sp{-1}\Gamma (1 + \nu d\sb 1 + s)\sp 3
\Gamma (1 + \nu d \sb 2 - s)\sp {-3}.
$$ 
Here $p = q = 6, \, m = 5, \, n = 3,$
$a\sb 1 = a\sb 2 = a\sb 3 = -\nu d\sb 1, \, 
a\sb 4 = a\sb 5 =  a\sb 6 = 1 + \nu d \sb 2,$ 
$b\sb 1 = b\sb 2 = b \sb 3 = 0, \, 
b\sb 4 = b\sb 5 = b\sb 6 = \nu, \, \Delta \sp * = -3\nu - 3;$
convergence conditions (C2) and (C3) are satisfied since 
now $\vert z \vert  \ge 1.$ The set of all the unremovable singular 
points of the function $g \sb {6, 6} \sp{(5, 3)} (s),$ encircled 
by the curve $L\sb 2,$ consists of the points 
$s  = -1 - \nu d\sb 1 - k$ with $k \in  \Bbb N - 1 ;$ 
each of these points is a pole of the second order. 
Therefore 
$$
\Res(g\sb {6, 6} \sp {(5, 3)}; - \nu d\sb{1} - 1 - k) =
\lim \limits \sb {s \to  - \nu d\sb{1} - 1 - k} 
\frac{\partial}{\partial s}\big((s + \nu d\sb{1} + 1 + k)\sp 2 
g\sb {6, 6} \sp {(5, 3)}\big),
$$ 
where $k \in \Bbb N - 1 .$ 

Let $s = - \nu d\sb 1  - 1 - k + u$ and 
$$H\sb 1(u) = 
g\sp {(5, 3)} \sb {6, 6}(-\nu d\sb 1  - 1 - k + u) = 
z\sp {- \nu d\sb 1 - 1 - k + u}\,\Gamma(\nu d\sb 1 + 1 + k - u)\sp 3 
\times 
$$
$$ 
\Gamma(\nu \Delta  + 1 + k - u)\sp 2\Gamma (-\nu \Delta - k + u)\sp{-1}
\Gamma (- k + u)\sp3 \Gamma(2\nu \Delta  + 2 + k - u)\sp {-3} =
$$ 
$$
= z\sp{- \nu d\sb 1 - 1 - k + u}\,\Gamma(\nu d\sb 1 + 1 + k - u)\sp 3 
\Gamma(\nu \Delta  + 1 + k - u)\sp3 \Gamma(\nu \Delta + 1 + k - u)\sp{-1} 
\times
$$
$$
\Gamma (-\nu \Delta - k + u)\sp{-1}\Gamma (1 + k - u)\sp3 
\Gamma (- k + u)\sp3 \Gamma (1 + k - u)\sp{-3} 
\Gamma (2 + 2\nu \Delta  + k - u)\sp{-3} = 
$$
$$
= \big(\pi/\sin(\pi u)\big)\sp{2} (-1)\sp {\nu \Delta}
z\sp{-\nu d\sb 1 - 1 - k + u}\, \Gamma(\nu d\sb 1 + 1 + k - u)\sp3 
\times
$$ 
$$
\Gamma(\nu \Delta + 1 + k - u)\sp3 \Gamma (1 + k - u)\sp{-3}
\Gamma (2 + 2\nu \Delta  + k - u)\sp{-3} =
\big(\pi/\sin(\pi u)\big)\sp{2} H\sp\ast(u),
$$ 
where 
$$
H\sp\ast(u) = (-1)\sp{\nu \Delta}\, z\sp{-\nu d\sb1 -1-k+u}\times
$$ 
$$
\Gamma (\nu d\sb 1  + 1 + k - u)\sp 3 \Gamma (1 + k - u)\sp {-3}
\Gamma (1 +  \nu \Delta  + k - u)\sp 3 \Gamma (2 + 2\nu \Delta \
+ k - u)\sp {-3}  =
$$ 
$$
= (-1)\sp {\nu \Delta } z\sp{- T + \nu}\Gamma(T)\sp 3 
\Gamma (T + 1 + \nu \Delta)\sp {-3}\Gamma (T - \nu \Delta)\sp{-3}
\Gamma (T - \nu )\sp 3  = 
$$ 
$$
= (-1)\sp {\nu \Delta} z\sp {- T + \nu} \bigg(\prod \limits \sb {\kappa = 
\nu  + 1} \sp {\nu \Delta} (T - \kappa)\bigg)\sp{\!3} 
\bigg(\prod\limits \sb {\kappa = 0} \sp {k + \nu \Delta}(T +\kappa)
\sp {-1}\bigg)\sp{\!3} =
$$ 
$$
= (-1)\sp {\nu \Delta} z\sp {- T + \nu} R\sb 0(T;\nu)
\sp 3 ((\nu \Delta )!/(\nu d\sb1)!)\sp{-3}
$$ 
and $T = \nu \Delta  + 1 + k - u.$ Therefore 
$$
(-1)\sp {\nu d\sb 2} ((\nu\Delta)!/(\nu d\sb1)!)\sp 3 
\Res(g\sp {(5, 3)} \sb {6, 6}; - 1 - \nu d\sb 1 - k ) = 
$$ 
$$
= (-1)\sp{\nu} z\sp {- T + \nu}\bigg(R\sb 0 (T;\nu)\sp 3 \log{z} 
- \frac{\partial}{\partial T} R\sb 0 (T;\nu )\sp 3\bigg) 
\bigg\vert \sb {T = 1 + \nu \Delta  + k} 
$$ 
because $(\pi u/(\sin(\pi u))\sp 2$ is an even function. Thus 
$$
f\sb 3 (z, \nu ) = f\sb 2 (z, \nu)\log{z} - (-1)\sp {\nu} ((\nu 
\Delta )!/(\nu d\sb1)!)\sp{-3}\cdot \hskip-4pt
\sum \limits \sb {t  = \nu \Delta  + 1} \sp \infty z\sp {-t + \nu} 
\frac{\partial}{\partial t}R\sb0(t;\nu)\sp3 ;
$$ 
since $R\sb 0 (t;\nu)\sp 3$ has zeros of the third order in the 
points $t = 1, \ldots, \nu \Delta,$ it follows that 
$$
f\sb 3 (z, \nu) = f\sb 2 (z, \nu)\log{z} - (-1)\sp \nu
((\nu \Delta)!/(\nu d\sb1)!)\sp{-3}\hskip-4pt\sum \limits \sb {t  = 1 + \nu} 
\sp {\infty} z\sp {-t + \nu}\frac{\partial}{\partial t}
R\sb 0(t;\nu)\sp3. \tag{2.9}
$$ 
Let
$$
f\sb 4 (z, \nu ) = - ((\nu \Delta )! / (\nu d\sb 1 )!)
\sp {-3}(-z)\sp {\nu}\hskip-4pt \sum \limits \sb {t = \nu \Delta  + 1} 
\sp \infty z\sp {-t}\frac{\partial}{\partial t}
R\sb0(t;\nu)\sp3; \tag{2.10}
$$ 
then 
$$
f\sb 3(z,\nu) = f\sb2 (z,\nu)\log{z} + f\sb 4 (z, \nu ). \tag{2.11}
$$ 
Let 
$$
f\sb 5 \sp {\vee} (z, \nu )  = -\frac{(-1)\sp{\nu \Delta}}{2\pi i} 
\int \limits \sb {L \sb 2} g \sb {6, 6} \sp{(6, 3)}(s)ds = 
$$ 
$$
= -(-1)\sp{\nu \Delta}\,G\sb{6,6}\sp{(6, 3)}
\left(- z \bigg\vert\matrix -\nu d\sb1, &\!\! -\nu d\sb1, &\!\! -\nu d\sb1, 
&\!\! 1+\nu d\sb2, &\!\! 1+\nu d\sb2, &\!\! 1+\nu d\sb2\\ 
0, &\!\! 0, &\!\! 0, &\!\! \nu, &\!\! \nu, &\!\! \nu\\
\endmatrix\right), 
$$ 
where $z \in \Omega \sp \ast \sb0, \, \nu \in \Bbb N,$ and 
$$
g \sb {6, 6} \sp{(6, 3)} = g \sb {6, 6} \sp{(6, 3)} (s) = 
(-z)\sp s \Gamma (-s)\sp 3 \Gamma (\nu  - s)\sp 3 
\Gamma (1  +  \nu d\sb 1 + s)\sp3 
\Gamma (1  +  \nu d\sb 2  -  s)\sp { - 3} .
$$ 
Here $p = q = 6, \, m = 6, \, n = 3,$ 
$a\sb 1 = a\sb 2 = a\sb 3 = -\nu d\sb 1, \, 
a\sb 4 = a\sb 5 =  a\sb 6 = 1 + \nu d \sb 2,$
$b\sb 1 = b\sb 2 = b \sb 3 = 0, \, 
b\sb 4 = b\sb 5 =  b\sb 6 = \nu, \, 
\Delta \sp * = -3\nu - 3;$ since $-z \in \Omega \sp \ast,$ each of 
convergence conditions (C2) and (C3) is satisfied. The set of all 
the unremovable singular points of the function $g \sb {6, 6} \sp{(5, 3)}(s),$ encircled by the curve $L\sb 2,$ consists of the points 
$s  = -1 - \nu d\sb 1 - k$ with $k \in  \Bbb N - 1 ;$ each of these 
points is a pole of the third order. Therefore
$$
\Res(g\sb {6, 6} \sp {(6, 3)}; - \nu d\sb{1} - 1 - k) = 
\lim \limits \sb {s \to  - \nu d\sb{1} - 1 - k} 
\frac12\big(\frac{\partial}{\partial s}\big)\sp 2
\big((s + \nu d\sb{1} + 1 + k)\sp 2 g\sb {6, 6} \sp {(5, 3)}\big),
$$ 
where $k \in \Bbb N - 1 .$ 

Let $s  = - \nu d\sb 1  - 1  -  k +  u,$ and 
$$
H\sb 2 (u)  = g\sp {(6, 3)} \sb {6, 6}(- \nu d\sb 1 - 1 - k + u)  =
$$
$$
= (-z)\sp {- \nu d\sb 1 - 1 - k + u}\Gamma(-k+ u)\sp 3
\Gamma (\nu d\sb 1  + 1  +  k - u)\sp 3 
\Gamma ( \nu \Delta + 1 + k - u)\sp 3 \Gamma(2  + 2\nu \Delta   
+  k  -  u)\sp {-3} 
$$
$$
= (-1)\sp k (-z)\sp {- \nu d\sb 1  - 1  - k  + u} 
\Gamma(1 + \nu \Delta   +  k - u)\sp 3 \Gamma (2 + 2\nu 
\Delta + k - u)\sp{-3}\times
$$ 
$$ 
\times\,\Gamma (\nu d\sb1  +  1  + k - u)\sp 3 
\Gamma (1  +  k - u)\sp{-3}(\pi /\sin(\pi u))\sp 3 
= (\pi /\sin(\pi u))\sp 3H\sb 2 \sp \ast (u),
$$ 
where 
$$
H\sb 2 \sp \ast (u) = (-1)\sp k(-z)\sp {-T + \nu}
R\sb 0 (T;\nu)\sp 3 ((\nu \Delta )!/(\nu d\sb1)!)\sp {-3} 
$$ 
and $T  = \nu \Delta   + 1  + k - u, \; t  = \nu \Delta  +  1  + k,$ 
as before. 

Therefore 
$$
-(-1)\sp {\nu \Delta}((\nu \Delta )!/(\nu d\sb1)!)\sp3 
\Res(g\sb {6, 6} \sp {(6, 3)} ; -  1 - \nu d\sb 1 - k)  =
$$ 
$$
= -(-1)\sp {\nu \Delta  + k} 
\bigg(\frac12((\log( - z))\sp 2 (-z)\sp {-t  + \nu}R\sb 0 (t;\nu)\sp 3  - 
(-z)\sp{-t+\nu}\log(-z)\frac{\partial}{\partial T}R\sb 0(t;\nu)\sp 3 + 
$$
$$
+ \frac12\big((-z)\sp{-t + \nu}\big(\frac{\partial}{\partial T}\big)\sp2 
R\sb0(t;\nu)\sp3 + \pi\sp2(-z)\sp{ - t  +  \nu} R\sb 0 (t;\nu)\sp3
\big)\bigg)
$$ 
because $(\pi u/\sin(\pi u))\sp3  =
1  + \frac12(\pi u)\sp 2 + \ldots.$ 

Consequently,  
$$
f\sb5\sp\vee (z, \nu) = \frac12f\sb 2(z,\nu)\big((\log(-z))\sp 2   
+  \pi \sp 2\big) + f\sb 4 (z, \nu)\log(-z) +  \tag{2.12}
$$ 
$$
+ \frac12(-z)\sp \nu ((\nu \Delta )!/(\nu d\sb 1)!)\sp {-3} 
\sum \limits \sp \infty \sb {t   =  \nu \Delta  + 1} z\sp{- t} 
\big(\frac{\partial}{\partial T}\big)\sp 2 R\sb 0(t;\nu)\sp 3  =
$$ 
$$
= \frac12\,f\sb 2 (z, \nu)(\log{z})\sp 2 + f\sb 4 (z, \nu)\log{z} + 
$$
$$
+ \frac12(-z)\sp \nu ((\nu \Delta )!/(\nu d\sb 1)!)\sp {-3} 
\hskip-4pt\sum \limits \sb {t = \nu \Delta  + 1} \sp \infty z\sp{-t} 
\big(\frac{\partial}{\partial T}\big)\sp2 R\sb 0(t;\nu)\sp 3 
- \pi if\sb 3 (z, \nu),
$$
where $z\in \Omega \sb 0 \sp *.$ 

Since $R\sb0(t;\nu)\sp 3$ has zeros of the 
third order at the points $t = 1 + \nu, \ldots, \nu \Delta,$ it follows 
that 
$$
f\sp\vee \sb 5 (z, \nu) = \frac12 f\sb2 (z,\nu) (\log{z})\sp2 + 
f\sb4(z,\nu)\log{z}\, +
$$ 
$$
+ \frac12(-z)\sp \nu ((\nu \Delta)!/(\nu d\sb1)!)\sp {-3} 
\sum \limits \sp\infty \sb {t = \nu  + 1}z\sp {-t}
\big(\frac{\partial}{\partial T}\big)\sp2 R\sb 0(t;\nu)\sp3 - 
i\pi f\sb 3(z, \nu).
$$ 
Let 
$$
f\sb6(z, \nu) = \frac12(-z)\sp \nu ((\nu \Delta)!/(\nu d\sb1)!)\sp{-3} 
\sum \limits \sp \infty \sb {t = 1} 
z\sp{-t}\big(\frac{\partial}{\partial T}\big)\sp2 R\sb 0(t;\nu)\sp 3, 
\tag{2.13} 
$$
$$
f\sb 5(z, \nu) = f\sb 6 (z, \nu) + 
\frac12 f\sb 2 (z,\nu)(\log{z})\sp2 + f\sb4 (z, \nu)\log{z},
$$ 
where $z \in \Omega \sp \ast \sb 0.$ Then 
$$
f\sb5 (z, \nu) = f\sp \vee \sb5(z, \nu) + i\pi f\sb3 (z, \nu), 
\tag{2.14}
$$ 
where $z \in \Omega \sp \ast \sb 0.$ Let further 
$$
f\sp \ast \sb j (z, \nu) = \big((\nu \Delta)!/(\nu d\sb1)!\big)\sp3  
f\sb j (z, \nu),  \qquad j = 1, \ldots, 6.
\tag{2.15}
$$
Expanding function $R\sb0(t;\nu)\sp3 $ into partial fractions, 
we obtain 
$$
R\sb 0(t;\nu)\sp3 = \sum \limits \sb{k = 0} \sp {\nu \Delta} 
\alpha \sp \ast \sb {\nu, k} (t + k) \sp {-3} + \sum \limits \sb{k = 0} 
\sp {\nu \Delta} \beta \sp \ast \sb {\nu, k} (t + k) \sp {-2} + 
\sum \limits \sb{k = 0} \sp {\nu \Delta} \gamma \sp \ast \sb {\nu, k} 
(t + k) \sp {-3} 
$$  
with 
$$
\alpha\sp\ast\sb{\nu, k} = (-1)\sp {\nu + \nu \Delta + k} 
\binom{\nu \Delta} {k}\sp{\!3}\binom {\nu \Delta + k}
{\nu \Delta  - \nu}\sp{\!3}, \tag{2.16}
$$  
$$
\beta \sp \ast \sb {\nu, k} = \lim \limits \sb {t \rightarrow 
 - k} \frac{\partial}{\partial t}\big(R\sb 0(t;\nu)\sp3 (t + k)\sp 3\big) = 
\alpha \sp \ast \sb {\nu, k}\big(3(-\hskip-8pt\sum \limits \sb {\kappa  = 
\nu  + k + 1} \sp {\nu \Delta + k} \hskip-5pt\kappa \sp {-1} - \sum \limits 
\sb {\kappa  = 1} \sp {\nu \Delta - k} \kappa \sp {-1} + \sum \limits 
\sb {\kappa  = 1} \sp k \kappa \sp {-1})\big), \tag{2.17} 
$$ 
$$
2\gamma \sp \ast \sb {\nu, k} = \alpha \sp \ast 
\sb {\nu, k}\big(3(-\hskip-8pt\sum \limits \sb {\kappa  = \nu  + k + 1} 
\sp {\nu \Delta + k} \hskip-5pt\kappa \sp {-1} - \sum \limits 
\sb{\kappa  = 1} \sp {\nu \Delta - k} \kappa \sp {-1} + 
\sum \limits \sb {\kappa  = 1}\sp k \kappa \sp {-1})\big)\sp 2 -  
\tag{2.18}
$$ 
$$
\alpha \sp \ast \sb {\nu, k}\big(3(\hskip-8pt\sum \limits \sb {\kappa  = \nu  
+ k + 1} \sp {\nu \Delta + k}\hskip-5pt \kappa \sp {-2} - \sum \limits 
\sb {\kappa  = 1} \sp {\nu \Delta - k} \kappa \sp {-2} - 
\sum \limits \sb {\kappa  = 1} \sp k \kappa \sp {-2})\big),
$$ 
where $k = 0, \ldots, \nu \Delta .$ 
It follows from (2.4), (2.8), (2.10), (2.13) and (2.15) - (2.18) that 
$$
f\sp \ast \sb 2(z, \nu) = (-z)\sp \nu \sum \limits 
\sb {t = 1 + \nu}\sp {+\infty}z\sp{-t}R\sb 0(t;\nu)\sp 3 = \tag{2.19}
$$ 
$$
= (-z)\sp \nu \sum \limits \sb {t = 1 + \nu} \sp {+\infty} z\sp {-t} 
\sum \limits \sb {k  =  0}\sp {\nu \Delta} \alpha \sp \ast \sb {\nu, k}
(t + k)\sp {-3} + (-z)\sp \nu \sum \limits \sb {t = 1 + \nu} 
\sp {+\infty} z\sp {-t} \sum \limits \sb {k  =  0}\sp {\nu \Delta} 
\beta \sp \ast \sb {\nu, k}(t + k)\sp {-2}\, + 
$$ 
$$
+\, (-z)\sp \nu \sum \limits \sb {t = 1 + \nu} \sp {+\infty} z\sp {-t} 
\sum \limits \sb {k  =  0}\sp {\nu \Delta} \gamma \sp \ast \sb {\nu, k}
(t + k)\sp {-1} = (-z)\sp \nu \sum \limits \sb {k  =  0}
\sp {\nu \Delta} \alpha \sp \ast \sb {\nu, k} z\sp k 
\sum \limits \sb {t = 1 + \nu  + k} \sp {+\infty} z\sp {-t} t
\sp {-3} \, + 
$$ 
$$
+\, (-z)\sp \nu \sum \limits \sb {k  =  0}\sp {\nu \Delta} \beta 
\sp \ast \sb {\nu, k} z\sp k \sum \limits \sb {t = 1 + \nu  + k} 
\sp {+\infty} z\sp {-t} t\sp {-2} + (-z)\sp \nu 
\sum \limits \sb {k  =  0}\sp {\nu \Delta} \gamma \sp \ast 
\sb {\nu, k} z\sp k \sum \limits \sb {t = 1 + \nu  + k} 
\sp {+\infty} z\sp {-t} t\sp {-1} = 
$$ 
$$
= \alpha \sp \ast (z;\nu)L\sb 3(z\sp {-1}) + \beta \sp \ast (z;\nu)L
\sb 2(z\sp {-1}) + \gamma \sp \ast (z;\nu) (-\log(1 - 1 / z)) - \phi 
\sp \ast(z;\nu) 
$$ 
with
\vskip-6pt
$$
L\sb n (z) = \sum \limits \sp {+ \infty} 
\sb {t = 1} z\sp t t\sp {-n}, \tag{2.20}
$$
\vskip-6pt
$$
\alpha \sp \ast (z;\nu) = (-z)\sp\nu \sum \limits 
\sb {k = 0} \sp {\nu \Delta} \alpha \sp\ast \sb {\nu, k}z\sp k = 
f\sp \ast \sb 1 (z;\nu), \tag{2.21}
$$ 
\vskip-6pt
$$
\beta \sp \ast (z;\nu) = (-z)\sp\nu \sum \limits \sb {k = 0} 
\sp {\nu \Delta} \beta \sp\ast \sb {\nu, k}z\sp k \tag{2.22}
$$ 
\vskip-6pt
$$
\gamma \sp \ast (z;\nu) = (-z)\sp\nu \sum \limits \sb {k = 0} 
\sp {\nu \Delta} \gamma \sp\ast \sb {\nu, k}z\sp k, \tag{2.23}
$$ 
\vskip-6pt
$$
\phi \sp\ast(z;\nu) = (-z)\sp \nu \sum \limits \sb {k = 0} 
\sp {\nu \Delta}\alpha \sp \ast \sb {\nu, k}z\sp k \sum \limits \sp {k 
+ \nu} \sb {t = 1}z\sp {- t}t\sp {-3}\, +  \tag{2.24}
$$ 
\vskip-6pt
$$
+\,(-z)\sp \nu \sum \limits \sb {k = 0} \sp {\nu \Delta}\beta 
\sp \ast \sb {\nu, k}z\sp k \sum \limits \sp {k + \nu} \sb {t = 1}z
\sp {- t}t\sp {-2} + (-z)\sp \nu \sum \limits \sb {k = 0} 
\sp {\nu \Delta}\gamma \sp \ast \sb {\nu, k}z\sp k \sum \limits 
\sp {k + \nu} \sb  {t = 1}z\sp {- t}t\sp {-1} =
$$ 
$$
= (-z)\sp \nu \sum \limits \sp \nu \sb {t = 1}z\sp{-t}
\big(t\sp {-3}\alpha \sp \ast (z;\nu) +t\sp {-2}\beta \sp \ast (z;\nu) 
+ t\sp {-1}\gamma \sp \ast (z;\nu)\big)\, +
$$ 
$$
+\,(-z)\sp \nu \sum \limits \sb {k = 0} \sp {\nu \Delta}\alpha 
\sp \ast \sb {\nu, k}z\sp k \sum \limits \sp {k + \nu} 
\sb {t = 1 + \nu}z\sp {- t}t\sp {-3} + (-z)\sp \nu 
\sum \limits \sb {k = 0} \sp {\nu \Delta}\beta \sp \ast 
\sb {\nu, k}z\sp k \sum \limits \sp {k + \nu} \sb {t = 1 + \nu}z
\sp {- t}t\sp {-2}\, +
$$ 
$$
+\,(-z)\sp \nu \sum \limits \sb {k = 0} \sp {\nu \Delta}\gamma \sp \ast 
\sb {\nu, k}z\sp k \sum \limits \sp {k + \nu} \sb {t = 1 + \nu}z
\sp {- t}t\sp {-1};
$$ 
\vskip-12pt
$$
f\sp \ast \sb 4(z, \nu) = -(-z)\sp \nu \sum \limits \sb 
{t = 1 + \nu}\sp {+\infty}z\sp{-t}\frac{\partial}{\partial t}
R\sb 0(t;\nu)\sp3 = \tag{2.25}
$$ 
\vskip-16pt
$$
= (-z)\sp \nu \sum \limits \sb {t = 1 + \nu} \sp {+\infty} z
\sp {-t} \sum \limits \sb {k  =  0}\sp {\nu \Delta} 3\alpha \sp \ast 
\sb {\nu, k}(t + k)\sp {-4} + (-z)\sp \nu \sum \limits \sb {t = 1 
+ \nu} \sp {+\infty} z\sp {-t} \sum \limits \sb {k  =  0}
\sp {\nu \Delta}2\beta \sp \ast \sb {\nu, k}(t + k)\sp {-3}\, + 
$$ 
\vskip-12pt
$$
+\,(-z)\sp \nu \sum \limits \sb {t = 1 + \nu} \sp {+\infty} z\sp {-t} 
\sum \limits \sb {k  =  0}\sp {\nu \Delta} \gamma \sp \ast 
\sb {\nu, k}(t + k)\sp {-2} = (-z)\sp \nu \sum \limits \sb {k  =  0}
\sp {\nu \Delta} 3\alpha \sp \ast \sb {\nu, k} z\sp k \sum \limits \sb 
{t = 1 + \nu  + k} \sp {+\infty} z\sp {-t} t\sp {-4}\, + 
$$ 
$$
+\,(-z)\sp \nu \sum \limits \sb {k  =  0}\sp {\nu \Delta} 2\beta 
\sp \ast \sb {\nu, k} z\sp k \sum \limits \sb {t = 1 + \nu  + k} 
\sp {+\infty} z\sp {-t} t\sp {-3} + (-z)\sp \nu \sum \limits \sb 
{k  =  0}\sp {\nu \Delta} \gamma \sp \ast \sb {\nu, k} z\sp k 
\sum \limits \sb {t = 1 + \nu  + k} \sp {+\infty} z\sp {-t} t
\sp {-2} = 
$$ 
$$
= 3\alpha \sp \ast (z;\nu)L\sb 4(z\sp {-1}) + 2\beta \sp \ast (z;\nu)
L\sb 3(z\sp {-1}) + \gamma \sp \ast (z;\nu) L\sb 2 (z)\sp {-1} - 
\psi \sp \ast(z;\nu), 
$$ 
where 
\vskip-16pt
$$
\psi \sp\ast (z;\nu) = (-z)\sp \nu \sum \limits \sb {k = 0} 
\sp {\nu \Delta}3\alpha \sp \ast \sb {\nu, k}\,z\sp k \sum \limits \sp {k 
+ \nu} \sb {t = 1}z\sp {- t}t\sp {-4}\, +  \tag{2.26}
$$
$$
+\,(-z)\sp \nu \sum \limits \sb {k = 0} \sp {\nu \Delta}2\beta 
\sp \ast \sb {\nu, k}\,z\sp k \sum \limits \sp {k + \nu} \sb {t = 1}z
\sp {- t}t\sp {-3} + (-z)\sp \nu \sum \limits \sb {k = 0} 
\sp {\nu \Delta}\gamma \sp \ast \sb {\nu, k}\,z\sp k \sum \limits 
\sp {k + \nu} \sb  {t = 1}\,z\sp {- t}t\sp {-2} =
$$ 
$$
= (-z)\sp \nu \sum \limits \sp \nu \sb {t = 1}z\sp{- t}
\big(t\sp {-4}3\alpha \sp \ast (z;\nu) +t\sp {-3}2\beta 
\sp \ast (z;\nu) + t\sp {-2}\gamma \sp \ast (z;\nu)\big)\, +
$$ 
$$
+\,(-z)\sp \nu \sum \limits \sb {k = 0} \sp {\nu \Delta}3\alpha 
\sp \ast \sb {\nu, k}z\sp k \sum \limits \sp {k + \nu} 
\sb {t = 1 + \nu}z\sp {- t}t\sp {-4} + (-z)\sp \nu 
\sum \limits \sb {k = 0} \sp {\nu \Delta}2\beta \sp \ast 
\sb {\nu, k}z\sp k \sum \limits \sp {k + \nu} \sb {t = 1 
+ \nu}z\sp {- t}t\sp {-3}\, +
$$ 
$$
+\,(-z)\sp \nu \sum \limits \sb {k = 0} \sp {\nu \Delta}\gamma 
\sp \ast \sb {\nu, k}z\sp k \sum \limits \sp {k + \nu} 
\sb {t = 1 + \nu}z\sp {- t}t\sp {-2};
$$ 
$$
f\sp \ast \sb 6(z, \nu) = \frac12(-z)\sp \nu \sum \limits 
\sb {t = 1 + \nu}\sp {+\infty}z\sp{-t}(\frac{\partial}{\partial t})
\sp2 R\sb 0(t;\nu)\sp 3 = \tag{2.27}
$$ 
$$
= (-z)\sp \nu \sum \limits \sb {t = 1 + \nu} \sp {+\infty} z
\sp {-t} \sum \limits \sb {k  =  0}\sp {\nu \Delta} 6\alpha 
\sp \ast \sb {\nu, k}(t + k)\sp {-5} + (-z)\sp \nu 
\sum \limits \sb {t = 1 + \nu} \sp {+\infty} z\sp {-t} 
\sum \limits \sb {k  =  0}\sp {\nu \Delta}3\beta \sp \ast 
\sb {\nu, k}(t + k)\sp {-4}\, +  
$$ 
$$
+\,(-z)\sp \nu \sum \limits \sb {t = 1 + \nu} \sp {+\infty} z
\sp {-t} \sum \limits \sb {k  =  0}\sp {\nu \Delta} \gamma \sp \ast 
\sb {\nu, k}(t + k)\sp {-3} = (-z)\sp \nu \sum \limits \sb {k  =  0}
\sp {\nu \Delta} 6\alpha \sp \ast \sb {\nu, k} z\sp k 
\sum \limits \sb {t = 1 + \nu  + k} \sp {+\infty} z\sp {-t} t
\sp {-5}\, + 
$$ 
$$
+\,(-z)\sp \nu \sum \limits \sb {k  =  0}\sp {\nu \Delta} 3\beta 
\sp \ast \sb {\nu, k} z\sp k \sum \limits \sb {t = 1 + \nu  + k} 
\sp {+\infty} z\sp {-t} t\sp {-4} + (-z)\sp \nu 
\sum \limits \sb {k  =  0}\sp {\nu \Delta} \gamma 
\sp \ast \sb {\nu, k} z\sp k \sum \limits \sb {t = 1 + \nu  + k} 
\sp {+\infty} z\sp {-t} t\sp {-3} = 
$$ 
$$
= 6\alpha \sp \ast (z;\nu)L\sb 5(z\sp {-1}) + 
3\beta \sp \ast (z;\nu)L\sb 4(z\sp {-1}) + \gamma \sp \ast 
(z;\nu) L\sb 3 (z\sp {-1}) - \xi \sp \ast(z;\nu),
$$  
and 
$$
\xi \sp\ast(z;\nu) = (-z)\sp \nu \sum \limits \sb {k = 0} 
\sp {\nu \Delta}6\alpha \sp \ast \sb {\nu, k}z\sp k \sum \limits \sp {k 
+ \nu} \sb {t = 1}z\sp {- t}t\sp {-5}\, + \tag{2.28}
$$ 
$$
+\,(-z)\sp \nu \sum \limits \sb {k = 0} \sp {\nu \Delta}3\beta 
\sp \ast \sb {\nu, k}z\sp k \sum \limits \sp {k + \nu} \sb {t = 1}z
\sp {- t}t\sp {-4} + (-z)\sp \nu \sum \limits \sb {k = 0} 
\sp {\nu \Delta}\gamma \sp \ast \sb {\nu, k}z\sp k \sum \limits 
\sp {k + \nu} \sb  {t = 1}z\sp {- t}t\sp {-3} =
$$ 
$$
= (-z)\sp \nu \sum \limits \sp \nu \sb {t = 1}z\sp {- t}
\big(6t\sp{-5}\alpha \sp \ast (z;\nu) + 3t\sp{-4}\beta\sp\ast(z;\nu) 
+ t\sp{-3}\gamma \sp \ast (z;\nu)\big)\, +
$$ 
$$
+\,(-z)\sp \nu \sum \limits \sb {k = 0} \sp {\nu \Delta}6\alpha 
\sp \ast \sb {\nu, k}z\sp k \sum \limits \sp {k + \nu} 
\sb {t = 1 + \nu}z\sp {- t}t\sp {-5} + (-z)\sp \nu 
\sum \limits \sb {k = 0} \sp {\nu \Delta}3\beta \sp \ast 
\sb {\nu, k}z\sp k \sum \limits \sp {k + \nu} \sb {t = 1 + \nu}z
\sp {- t}t\sp {-4}\, +
$$ 
$$
+\,(-z)\sp \nu \sum \limits \sb {k = 0} \sp {\nu \Delta}\gamma 
\sp \ast \sb {\nu, k}z\sp k \sum \limits \sp {k + \nu} 
\sb {t = 1 + \nu}z\sp {- t}t\sp {-3}.
$$ 

\vskip2cm

\address{ {\it E-mail:}{\sl\ gutnik{\@}gutnik.mccme.ru}
}
\endaddress 
\enddocument